\documentclass[12pt]{amsart}

\usepackage[leqno]{amsmath}
\usepackage{amssymb}
\usepackage{amsxtra}
\usepackage{amscd}
\usepackage{hyperref}
\usepackage{graphicx, color}
\usepackage{dcpic}

\newtheorem{thm}{Theorem}
\newtheorem{prop}[thm]{Proposition}

\newtheorem{defn}[thm]{Definition}

\newcommand{\R}{\mathbb{R}}

\newcommand{\Z}{\mathbb{Z}}

\newcommand{\mir}[1]{\overline{#1}}

\newcommand{\graph}[1]{\Gamma_{L}}
\newcommand{\circles}[1]{\ensuremath{\mathrm{\textsc{Circles}}(#1)}}
\newcommand{\cutcircles}[1]{\ensuremath{\mathrm{\textsc{CutCircs}}(#1)}}
\newcommand{\freecircles}[1]{\ensuremath{\mathrm{\textsc{FreeCircs}}(#1)}}

\newcommand{\lra}{\longrightarrow}
\newcommand{\inlinediag}[2][0.33]{\includegraphics[scale=#1]{./#2}}
\newcommand{\myfig}[3][0.5]{\begin{center}\begin{figure}\includegraphics[scale=#1]{./#2}\caption{#3}\label{fig:#2}\end{figure}\end{center}}

\newcommand{\cross}[1]{\ensuremath{\mathrm{\textsc{cr}}(#1)}}

\newcommand{\resolution}[1]{\ensuremath{\mathrm{\textsc{res}}(#1)}}

\newcommand{\lefty}[1]{\overleftarrow{#1}}
\newcommand{\righty}[1]{\overrightarrow{#1}}

\newcommand{\cleaved}[1]{\mathcal{C\!L}_{#1}}

\newcommand{\leftComplex}[1]{\langle\!\!\!\langle\,#1\,]\!]} 
\newcommand{\rightComplex}[1]{[\![\, #1 \,\rangle\!\!\!\rangle}
\newcommand{\complex}[1]{\ensuremath{\langle\!\!\!\langle\,#1\,\rangle\!\!\!\rangle}}

\newcommand{\Oz}{Ozsv\'{a}th}

\title{The decategorification of bordered Khovanov homology}
\author{Lawrence P. Roberts}

\begin{document}
\begin{abstract}    
In \cite{RobD}, \cite{RobA} the author showed how to decompose the Khovanov homology of a link $\mathcal{L}$ into the algebraic pairing of a type D structure and a type A structure (as defined in bordered Floer homology), whenever a diagram for $\mathcal{L}$ is decomposed into the union of two tangles. Since Khovanov homology is the categorification of a version of the Jones polynomial, it is natural to ask what the type A and type D structures categorify, and how their pairing is encoded in the decategorifications. In this paper, the author constructs the decategorifications of these two structures, in a manner similar to Ina Petkova's decategorification of bordered Floer homology, \cite{Petk}, and shows how they recover the Jones polynomial. We also give a new proof of the mutation invariance of the Jones polynomial which uses these decomposition techniques. 
\end{abstract}
\maketitle

\section{Background and motivation}

\noindent In \cite{Khov}, M. Khovanov describes a homology theory whose ``Euler characteristic'' is a reparametrization of the Jones polynomial of an oriented link $\mathcal{L}$ in $S^{3}$. In particular, he uses a link diagram $L$ for $\mathcal{L}$ and defines a bigraded, free Abelian group $\complex{L;\Z}^{\ast,\ast}$ whose generators are decorated resolutions of $L$. The group $\complex{L;\Z}^{\ast,\ast}$ admits a $(1,0)$ differential $\partial_{\textsc{Kh}}$ for which the homology $\textsc{Kh}^{\ast,\ast}(L,\Z)$ is an invariant of the link $\mathcal{L}$. \\
\ \\
\noindent Since the differential is $(1,0)$ we can decompose $\complex{L;\Z}^{\ast,\ast}$ as a direct sum of chain complexes $\bigoplus_{j \in \Z} (C^{\ast,j},\partial_{\textsc{KH}})$  where $C^{\ast,j}$ is the subgroup of $\complex{L;\Z}^{\ast,\ast}$ whose second grading, called the quantum grading, is equal to $j$. The ``Euler characteristic'' of $\complex{L;\Z}^{\ast,\ast}$ is taken to be the polynomial
$$
J_{L}(q) = \sum_{j \in \Z} \chi(C^{\ast,j}) q^{j}
$$
where $\chi(C^{\ast,j})$ is the usual Euler characteristic of a finitely generated chain complex with free chain groups. This polynomial has the form $(-1)^{n_{-}(L)}q^{(n_{+}-2n_{-})(L)}\widetilde{J}_{L}(q)$, where $n_{\pm}(L)$ is the number of positive/negative crossings in $L$ and $\widetilde{J}_{L}(q)$ can be computed from the unoriented link diagram.  If we let $U$ be an unknot, and $L \sqcup U$ be a link with an unlinked and unknotted component $U$, then $\widetilde{J}_{L}(q)$ satisfies
$$
\begin{array}{l}
\widetilde{J}_{U}(q) = q + q^{-1} \\
\widetilde{J}_{L \sqcup U}(q) = (q + q^{-1}) \widetilde{J}_{L} \\
\widetilde{J}_{L} = \widetilde{J}_{L_{0}} - q \widetilde{J}_{L_{1}} \\
\end{array}
$$
Thus $J_{L}(q)$ is a version of the unnormalized Jones polynomial. \\
\ \\
\noindent Due to this relationship, made more precise below, Khovanov homology is called a {\em categorification} of the Jones polynomial, and the polynomial is called the {\em decategorification} of Khovanov homology. These terms are used since the polynomial arises from a process more general than just taking Euler characteristics. Despite both the homology and its decategorification giving link invariants, the homology $\complex{L}$ is known to be stronger than the Jones polynomial, and has many additional properties related to the cobordism on links. \\
\ \\
\noindent In \cite{RobD}, \cite{RobA} the author, inspired by bordered Floer homology, describes a similar construction for tangles\footnote{There are several other constructions of Khovanov homology for tangles, \cite{Khta}, \cite{APS1}, \cite{APS}, \cite{BarN}, \cite{LaPf}}.  The formal algebraic structure of these tangle invariants mimics the formal structure of \Oz, Lipshitz, and Thurston's description of bordered Floer homology, \cite{Bor1}, which also provided a road map for constructing a gluing theory for Khovanov homology.  \\
\ \\
\noindent The construction in \cite{RobA} takes a tangle diagram $T$ in a disc $D$, with a marked point $\ast \in \partial D$:\\

\begin{center}
\includegraphics[scale=1]{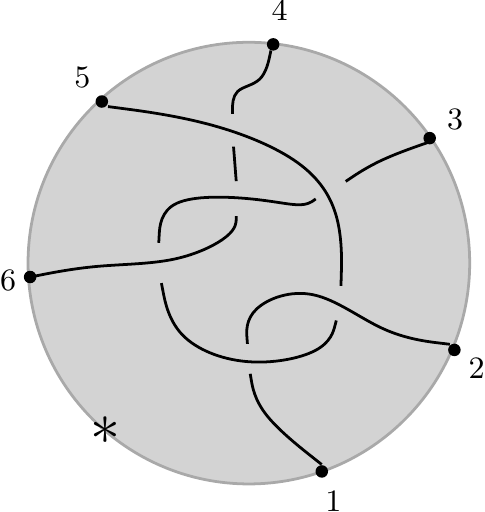}
\end{center}

\noindent and associates to it a bigraded differential module $\leftComplex{T}$ over a differential bigraded algebra $\mathcal{B}\Gamma_{n}$, where $T$ has $2n$ endpoints on $\partial D$. Each differential is $(1,0)$ in the respective bigradings, and the action $\leftComplex{T} \otimes_{I} \mathcal{B}\Gamma_{n} \lra  \leftComplex{T}$ preserves the bigradings on its source and target. The homotopy type of $\leftComplex{T}$ as an $A_{\infty}$-module is invariant under the three Reidemeister moves applied to swatches in the diagram $T$, and thus defines a tangle invariant. \\
\ \\
\noindent The main goal in defining these algebraic structures is to obtain a complete gluing theory for Khovanov homology. For instance, let $L$ be a link diagram in an oriented sphere $S$, and $S = D_{1} \cup D_{2}$, where each closed disk $D_{i}$ inherits its orientation from $S$ and $D_{1} \cap D_{2} = C$ is a circle in $S$. We will orient $C$ as the boundary of $D_{1}$. Suppose $L$ is transverse to $C$ away from its crossings. Then $L \cap D_{1} = T_{1}$ and $L \cap D_{2} = T_{2}$ are tangles in oriented disks. If we choose $\ast \in C$, then each disk also has a marked point in its boundary. We illustrate this situation in the diagram below, where we have taken $\ast$ to be the point at infinity. $D_{1}$ is then the left half plane, and $D_{2}$ is the right half plane.

\begin{center}
\includegraphics[scale=0.5]{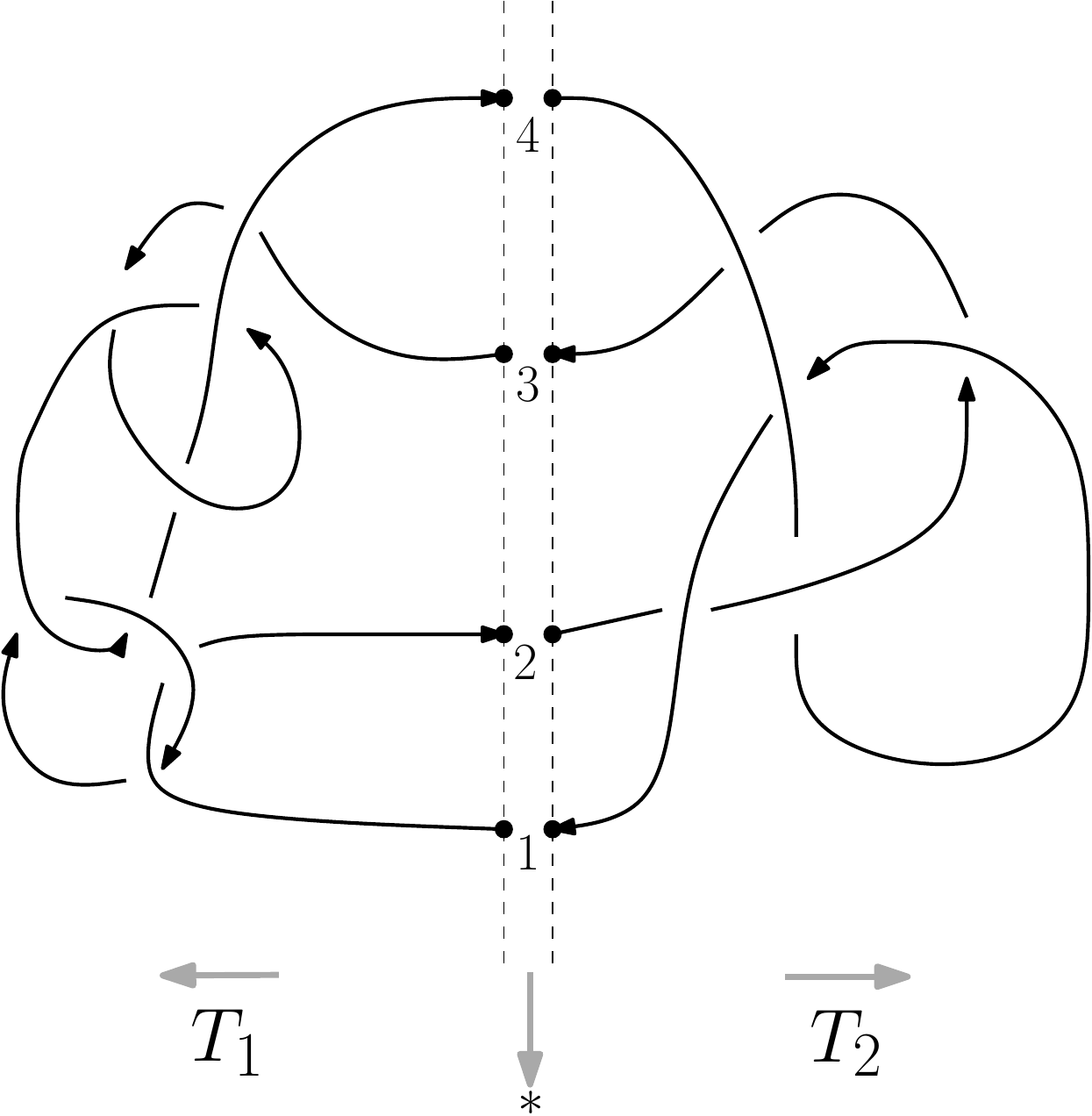}
\end{center}

\noindent To $(D_{1}, T_{1})$ assign the bigraded differential module $\leftComplex{T_{1}}$ described above. To $(D_{2}, T_{2})$ the construction in \cite{RobD} assigns a bigraded Abelian group $\rightComplex{T_{2}}$ and a $(1,0)$ map $\delta: \rightComplex{T_{2}} \lra  \mathcal{B}\Gamma_{n} \otimes_{I}  \rightComplex{T}$, which satisfies the structure relation for a type D structure, \cite{Bor1}. \\
\ \\
\noindent There is an algebraic construction $\boxtimes$ which pairs differential graded modules and type $D$ structures to obtain a normal chain complex. This can be adapted to the bigraded setting and applied to compute $\leftComplex{T_{1}} \boxtimes \rightComplex{T_{2}}$. The main result of \cite{RobA} is that $\leftComplex{T_{1}} \boxtimes \rightComplex{T_{2}}$ is chain homotopy equivalent to $\complex{L}$, the original link diagram. Furthermore, $\boxtimes$ preserves chain homotopy equivalence when we alter either factor by a homotopy equivalence (in the respective categories of $A_{\infty}$-modules and type $D$ structures), so this effects a gluing of the {\em homotopy types} of the invariants assigned to the tangles. As it is the homotopy types that are invariants of the tangles and links, and not the raw algebraic object, this provides a complete gluing theory for Khovanov homology.\\
\ \\
\noindent Having defined $\leftComplex{T_{1}}$ and $\rightComplex{T_{2}}$ it is natural to ask if there is any analog of these objects, and their pairing, in the simpler world of the Jones polynomial. In this paper we provide the answer to this questions by  describing the decategorifications of $\leftComplex{T}$ and $\rightComplex{T}$. The next section elaborates on the notion of decategorification we will use. In section \ref{sec:decat} the decategorification is described in the abstract. The section \ref{sec:idem} we provide more detail about the algebras involved, which allows us to give a concrete description of the decategorification in section \ref{sec:tangle}. This will all be described for the type A structure. In section \ref{sec:typeD} we show how to interpret these results for type D structures. Then we can show how to recover the Jones polynomial in section \ref{sec:Jones} and reprove its invariance under mutation in section \ref{sec:mutation}.     \\
\ \\   
\noindent The concrete description can be used to generalize this paper to other settings, and exhibit a kind of planar algebra structure, a topic the author pursues in a sequel, \cite{Plan}.\\
\ \\
\noindent{\em Convention:} If $M$ is a bigraded module, then $M\{(m,n)\}$ is the bigraded module with $(M\{(m,n)\})_{i,j} = M_{i-m,j-n}$. \\

\section{Background on decategorification}

\noindent First, we provide a little more detail concerning the method. Above, we referred to the polynomial $J_{L}(q)$ as a ``Euler characteristic.'' A better, and more sophisticated version of this statement can be obtained through the use of Grothendieck groups. For this paper, we will need a generalization of the following version of the Grothendieck construction, taken from \cite{Kho2}:

\begin{defn}
Let $A$ be a bigraded associative algebra. The Grothendieck  group $K_{0}(A)$ is the $\Z[q,q^{-1}]$-module generated by the elements $[P]$ where $P$ is a finitely generated, bigraded, projective $A$-module, and subject to the relations that $[P] = [P'] + [P'']$ when there is a short exact sequence $0 \lra P' \lra P \lra P'' \lra 0$, and $[P\{(m,n)\}]=(-1)^{m}q^{n}[P]$. 
\end{defn}

\noindent When $A$ is $\Z$ (in bigrading $(0,0)$)  the relations above confirm that $[\complex{L;\Z}^{\ast,\ast}] = \sum_{j} [C^{\ast,j}] = \sum_{j} q^{j} [C^{\ast,j}\{(0,-j)\}]$. Now $C^{\ast,j}\{(0,-j)\}$ is in gradings $(m,0)$ for $m \in \Z$. Thus, $[C^{\ast,j}\{(0,-j)\}] = \sum_{i} (-1)^{i}[C^{i,j}\{(-i,-j)\}]$. Since $\complex{L;\Z}^{\ast,\ast}$ is free, so is $C^{i,j}$ so the image of $C^{i,j}\{(-i,-j)\}$ (in bigrading $(0,0)$) in the Grothendieck group is just the rank of $C^{i,j}$. Thus,
$$
[\complex{L;\Z}^{\ast,\ast}] = \sum_{i,j \in Z} (-1)^{i}q^{j}\mathrm{rk}\,C^{i,j}
$$
from which it readily follows that $[\complex{L;\Z}^{\ast,\ast}] = J_{L}(q)$. \\
\ \\
\noindent This can be extended to the case where $P$ is a finitely generated, bigraded projective $A$-module with a $(1,0)$-differential. Following the pattern, the decategorifications of $\leftComplex{T}$ and $\rightComplex{T}$ should be the elements in a Grothendieck group for the differential bigraded algebra $\mathcal{B}\Gamma_{n}$ where $n$ is the number of unclosed components in $T$.\\
\ \\
\noindent Let $A$ be a bigraded differential graded algebra with $(1,0)$ differential. We will consider {\em right} differential graded modules over $A$, since that is the structure of $\leftComplex{T}$ over $\mathcal{B}\Gamma_{n}$. Following \cite{Kho2}, we aim to define the Grothendieck group $K_{0}(A)$. We need some additional definitions.

\begin{defn}
$\mathcal{K}(A)$ is the triangulated category found from the category of (right) bigraded differential modules with $(1,0)$ differential by quotienting out by the morphisms homotopic to zero. $\mathcal{D}(A)$ is the derived category found by localizing $\mathcal{K}(A)$ at its quasi-isomorphisms. 
\end{defn}

\noindent Here it is understood that the homotopies, chain maps, etc. do not change the quantum (second) grading. The distinguished triangles in the triangulated structure are those diagrams triangle isomorphic to a standard triangle
$$
M \stackrel{u}{\lra} N \lra C(u) \lra M\{(-1,0)\}
$$
where $u$ is a morphism of differential graded modules, and $C(u)$ is the mapping cone of $u$: the module $N \oplus M\{(-1,0)\}$ equipped with the differential 

$$
\left[\begin{array}{cc} d_{N} & u \\ 0 & -d_{M} \end{array}\right]
$$ 

\begin{defn}
A differential graded module $P$ is said to be {\em projective} if, given any (right) differential graded module $M$ with trivial homology, the complex $\textsc{Hom}(P,M)$, defined in \cite{Beil}, has trivial homology.     
\end{defn}

\begin{defn}
$\mathcal{K}\mathcal{P}(A)$ is the full subcategory of $\mathcal{K}(A)$ whose objects are the projective differential graded modules. 
\end{defn}

\noindent In Part II, section 10 of \cite{Beil} a bar construction is described which takes any differential graded module $M$ and finds a quasi-isomorphic projective differential graded module $B(M)$. Thus, $\mathcal{K}\mathcal{P}(A)$ is equivalent to the derived category $\mathcal{D}(A)$. 
 
\begin{defn}
An object $M$ of $\mathcal{K}\mathcal{P}(A)$ is {\em compact} if the natural inclusion
$$\bigoplus_{i \in I} \textsc{Hom}(M,N_{i}) \hookrightarrow \textsc{Hom}(M, \bigoplus_{i\in I} N_{i})$$ is an isomorphism for any collection of objects $\big\{N_{i}\,\big|\, i \in I\,\big\}$. 
\end{defn}

\begin{defn}
$\mathcal{P}(A)$ is the full subcategory of $\mathcal{K}(A)$ whose objects are compact, projective modules over $A$
\end{defn}

\noindent We are now in a position to define $K_{0}(A)$.

\begin{defn}
Let $A$ be a bigraded differential graded algebra with $(1,0)$ differential. $K_{0}(A)$ is the Grothendieck group of the category $\mathcal{P}(A)$. More specifically, $K_{0}(A)$ is the Abelian group with a generator $[P]$ for each compact, projective differential bigraded \rm{(}right\rm{)} module $P$ over $A$, with $(1,0)$ differential, subject to the relations $[P_{2}] = [P_{1}] + [P_{3}]$
for each distinguished triangle $P_{1} \rightarrow P_{2} \rightarrow P_{3}$
\end{defn}

\noindent It follows from the definition that $[P\{(1,0)\}] = -[P]$ (from the distinguished triangle coming from the mapping cone of the identity on $P$), and that $K_{0}(A)$ is a $\Z[q,q^{-1}]$-module, where the action is $q^{i}[P] = [P\{(0,i)\}]$. In particular, $[P\{(m,n)\}] = (-1)^m q^{n} [P]$. \\

\section{The Grothendieck group of $\mathcal{B}\Gamma_{n}$}\label{sec:decat}

\noindent We apply this construction to $\mathcal{B}\Gamma_{n}$. More details about this algebra will be given below. For this computation all that is required are the following properties, \cite{RobD}:
\begin{enumerate}
\item $\mathcal{B}\Gamma_{n}$ is the quotient of a quiver algebra $\mathcal{Q}\Gamma_{n}$ defined by an {\em acyclic} directed graph $\Gamma_{n}$,
\item The relations defining this quotient consist of identities involving paths of length $\geq 1$,
\item The $(1,0)$ differential is non-trivial only on paths of length $\geq 1$.
\end{enumerate} 

\noindent Let $I_{2n}$ be the algebra of (orthogonal) idempotents in $\mathcal{B}\Gamma_{n}$. We will now prove that

\begin{thm}\label{thm:KO}
$K_{0}(\mathcal{B}\Gamma_{n})$ is isomorphic to the $\Z[q^{1/2}, q^{-1/2}]$ module spanned by the idempotents corresponding to vertices in $\Gamma_{n}$.
\end{thm}

\noindent Note that this is in keeping with the computation of Grothendieck groups for acyclic quiver algebras with relations, \cite{Craw}. We have switched to $\Z[q^{1/2}, q^{-1/2}]$ since the quantum grading on $\mathcal{B}\Gamma_{n}$ is half-integral, but no other changes are necessary to the above construction.\\
\ \\
\noindent{\bf Proof:} Let $P$ be a right bigraded differential module over $\mathcal{B}\Gamma_{n}$. Since $P$ is projective and compact in $\mathcal{K}(\mathcal{B}\Gamma_{n})$ we may use a homotopy equivalent representative of $P$ which is finitely generated. Since $\Gamma_{n}$ is acyclic and directed, there is a vertex $v$ which has no out edges. Let $I_{v}$ be the corresponding idempotent. now consider $PI_{v}$. This is a submodule of $P$ since the action of any element of $\mathcal{B}\Gamma_{n}$ arises as the image of the action of path elements in $\mathcal{Q}\Gamma_{n}$. As no path starts at $v$, the only element which acts non-trivially on $PI_{v}$ is $I_{v}$. Furthermore, $PI_{v}$ is a subcomplex of $P$ with its $(1,0)$ differential since
$$
d_{P}(p \cdot I_{v}) = d_{P}(p)\cdot I_{v} + (-1)^{\ast}p \cdot d_{\mathcal{B}\Gamma_{n}} I_{v} =  d_{P}(p)\cdot I_{v}
$$
for any $p \in P$. Thus the image of any element in $PI_{v}$ under $d_{P}$ will be in $PI_{v}$. Thus there is a distinguished triangle
$$
P/PI_{v}\{(1,0)\} \lra PI_{v} \lra P \lra P/PI_{v}
$$
in $\mathcal{K}(\mathcal{B}\Gamma_{n})$. Therefore, $[PI_{v}]=[P] + [P/PI_{v}\{(1,0)\}] = [P] - [P/PI_{v}]$, so $[P] = [PI_{v}] + [P/PI_{v}]$. Now $P/PI_{v}$ is still a module over the differential graded algebra $B\Gamma_{n}$, but only the elements of $\mathcal{Q}\Gamma'$ will act non-trivially, where $\Gamma'$ is the directed graph $\Gamma_{n}\backslash\{v\}$. Since $\Gamma'$ is acyclic, it also has a vertex $v'$ with no outward edges. Due to the orthogonality of the idempotents $(P/PI_{v}) \cdot I_{v'} = PI_{v'}$. Applying this reasoning repeatedly, and using that $P$ is finitely generated, we arrive at
$$
[P] = \sum_{v \in \mathrm{vert}(\Gamma_{n})} [PI_{v}]
$$
Since the action of $\mathcal{B}\Gamma_{n}$ is essentially trivial on $PI_{v}$, we can consider $PI_{v}$ to be a bigraded complex over $\Z$ with $(1,0)$ differential, just as above. Thus $[PI_{v}] \in \Z[q^{1/2},q^{1/2}] [P_{v}]$ where $P_{v}$ is the module with a $\Z$ in bigrading $(0,0)$, trivial differential, and $P_{v} I_{v} = P_{v}$. Thus we see that the $\Z[q^{1/2}, q^{-1/2}]$ module spanned by the idempotents of $\mathcal{B}\Gamma_{n}$ maps surjectively to $K_{0}(\mathcal{B}\Gamma_{n})$. \\
\ \\
\noindent Before continuing we note that this is precisely what happens for bound quiver algebras, \cite{Craw}, as they too have Jordan-H\"older sequences of this type. \\
\ \\
\noindent We show that this map is an isomorphism by constructing a module $P$ for which $[P] =  \sum_{v \in \mathrm{vert}(\Gamma_{n})} J_{v}(q)[P_{v}]$ for any collection of polynomials $J_{v}(q)$ with integer coefficients. For each term $\pm a q^{j}$ we have a copy of $P_{v}^{a}\{(0,j)\}$ if the sign is $+$, and $P_{v}^{a}\{(1,j)\}$ if the sign is $-$. We take the direct sum over all such terms. The differential is taken to be trivial. We do this for all $v \in \mathrm{vert}(\Gamma_{n})$, and then define the action of $\mathcal{B}\Gamma_{n}$ to be trivial for every element that is the image of a path of $\mathcal{Q}\Gamma_{n}$ of length $\geq 1$, and let $I_{v}$ act nontrivially only on the copies of $P_{v}$. It is clear that this cannot be simplified further, and has   image $\sum_{v \in \mathrm{vert}(\Gamma_{n})} J_{v}(q)[P_{v}]$ in $K_{0}(\mathcal{B}\Gamma_{n})$. Thus there are no relations among the elements $[P_{v}]$. $\Diamond$.\\
\ \\
\noindent This provides a strategy for computing the image $[\leftComplex{T}]$, which we will describe presently. First, we spend some times on the idempotent sub-algebra $I_{2n}$. $\Diamond$
    
\section{The Idempotent sub-algebra of $\mathcal{B}\Gamma_{n}$}\label{sec:idem}

\noindent The vertices of the quiver defining $\mathcal{B}\Gamma_{n}$ correspond to certain planar configurations of circles and decorations, called {\em cleaved links} in \cite{RobD} and \cite{RobA}. More specifically, let $S$ be an oriented two-dimensional sphere. \\

\begin{defn}
A cleaved link $L$ in $S$ consists of the following data
\begin{enumerate}
\item a smoothly embedded circle $E_{L} \subset S$, called the {\em equator} for $L$,
\item a marked point $\ast_{L} \in E_{L}$
\item an identification of the closures of the two open discs $S \backslash E_{L}$ as $\lefty{D}_{L}$ and $\righty{D}_{L}$, called the inside and outside discs, respectively, where each is oriented from $S$, 
\item the orientation of $E_{L}$ induced as the boundary of $\lefty{D}_{L}$, and
\item a (possibly empty) set $\cutcircles{L}$ of $k_{L} \geq 0$ simple closed curves which each non-trivially and transversely intersect $E_{L}$ away from $\ast_{L}$
\end{enumerate}
\end{defn} 

\begin{defn}
For $L$ a cleaved link in an oriented sphere $S$, $P_{L}$ is the ordered set whose elements are points of intersection between $E_{L}$ and $\cutcircles{L}$ equipped with the ordering inherited from the orientation of $E_{L} \backslash \{\ast_{L}\}$. The cardinality of $P_{L}$ is $2n_{L}$.
\end{defn}
\ \\
\noindent Note that $P_{L}$ will be ordered opposite the orientation of $\righty{D}_{L}$. \\
\ \\
\noindent A cleaved link $L_{1}$ in a sphere $S_{1}$ is equivalent to another cleaved link $L_{2}$ in $S_{2}$ if there is an orientation preserving diffeomorphism $\phi: S_{1} \lra S_{2}$ which preserves each of the structures in the definition. In particular,\\

\begin{enumerate}
\item $\phi(E_{L_{1}}) = E_{L_{2}}$ and $\phi(\ast_{L_{1}})= \ast_{L_{2}}$,
\item $\phi(\lefty{D}_{L_{1}}) = \lefty{D}_{L_{2}}$, and thus $\phi(\righty{D_{L_{1}}}) = \righty{D_{L_{2}}}$
\item $\phi$ maps each $C \in \cutcircles{L_{1}}$ diffeomorphically to a circle $C' \in \cutcircles{L_{2}}$
\end{enumerate}
\ \\
\noindent It follows that $\phi$ induces an order preserving bijection $P_{L_{1}} \leftrightarrow P_{L_{2}}$.\\

\begin{defn}
A decorated, cleaved link $(L,\sigma)$ is a cleaved link $L$ and a map $\sigma : \cutcircles{L} \lra \{+,-\}$. The map $\sigma$ is called the {\em decoration}.
\end{defn}
\ \\
\noindent Two decorated, cleaved links are equivalent if there is an equivalence $\phi$ of the undecorated cleaved links with $\phi^{\ast}(\sigma_{2}) = \sigma_{1}$. \\

\begin{defn}
The set of equivalence classes of decorated, cleaved links $(L,\sigma)$ with $n_{L} = n$ will be denoted $\cleaved{n}$.
\end{defn}

\noindent In $\mathcal{B}\Gamma_{n}$ there is an idempotent $I_{(L,\sigma)}$ for each equivalence class of decorated, cleaved links. By theorem \ref{thm:KO}, we know that we should be interested in the $\Z[q^{1/2}, q^{-1/2}]$ modules spanned by these equivalence classes.\\

\begin{defn}
For each $n \geq 0$, $\mathcal{I}_{2n}$ is the free $\Z[q^{1/2}, q^{-1/2}]$-module generated by the elements of $\cleaved{n}$. The generator corresponding to $(L,\sigma) \in \cleaved{n}$ will be denoted $I_{(L,\sigma)}$. 
\end{defn}
\ \\
\noindent{\bf Examples:} When $n=0$, $\mathcal{I}_{0} = \Z[q^{1/2}, q^{-1/2}]$. It is possible to include this in the framework above, by allowing $\cutcircles{L}$ to be empty. Then $\mathcal{I}_{0}$ has a generator $I_{0}$ corresponding to the equivalence class for $S^{2} \subset \R^3$, oriented as the boundary of the unit ball, with $E_{L} = \{(x,y,0)|x^{2} + y^{2} = 1\}$, and $\ast = (1,0,0)$. We take $\lefty{D}$ to be the upper hemisphere, since that endows $E_{L}$ with the same orientation it inherits from being the boundary of $D^{2}$ in $\R^{2} \times \{0\} \subset \R^{3}$. However, we do not obtain a different cleaved link by taking $\lefty{D}$ to be the lower hemisphere since $\phi(x,y,z) = (x, -y, -z)$ is orientation preserving when restricted to $S^{2}$, takes $(1,0,0)$ and $E$ to themselves, and carries the upper hemisphere to the lower hemisphere.\\
\ \\ 
\noindent{\em Convention for describing generators:} Before giving more examples we describe how the choice of $\ast$, and the ordering of $P_{L}$, allows each generator to be identified by combinatorial data. Each generator is determined by two planar matchings of $P_{L}$: a planar matching $\lefty{m}_{L}$ embedded in $\lefty{D}_{L}$ and the other $\righty{m}_{L}$ embedded in $\righty{D}_{L}$. A planar matching on $2n_{L}$ enumerated points is uniquely determined by a permutation of the even numbers $2,4,\ldots, 2n_{L}$ by the rule that the $k^{th}$ even number in the permutation $2\sigma_{1}, \ldots, 2\sigma_{n_{L}}$ is the endpoint $p_{2\sigma_{k}}$ of the arc starting at $p_{2k-1}$. We will describe both $\lefty{m}_{L}$ and $\righty{m}_{L}$ by these permutations, as specified by $P_{L}$. To finish encoding $L$ we need to specify the decoration $\pm$ on each circle in $\cutcircles{L}$. We do this by first ordering the circles by the order in which we first meet the circles if we start at $\ast$ and walk around $E_{L}$ according to its orientation. Thus the circle containing $p_{1}$ will always come first in our ordering. This is equivalent to the rule $C < C'$ if, and only if, the smallest subscript of any $p_{j} \in P_{L}$ occurring in $C$ is less than the smallest subscript of any $p_{r}$ occurring in $C'$. With this ordering, a list of $k_{L}$ elements from $\{+,-\}$ corresponds to a choice of decorations on the circles of $\cutcircles{L}$:  the $i^{th}$ entry in the list is the decoration on the $i^{th}$ circle in the ordering. For example, in Figure \ref{fig:I4Gens} the generator $B_{-+}$ is specified by $\lefty{m}_{L} = (42)$ and $\righty{m}_{L} = (42)$. The $-$ decoration occurs on the circle through the point $1$, which is the larger circle in the picture, while the smaller, and second circle, is decorated with a $+$. \\
\ \\
\noindent The generators of $\mathcal{I}_{2}$ are cleaved links whose circles intersect its equator exactly twice. For each circle $C \subset S^{2}$ that intersects another circle $E$ there must be at least two intersections. Consequently $\cutcircles{L} = \{C\}$. $C$ can be decorated with either a $+$ or a $-$. Thus $\cleaved{2} = \{C_{+}, C_{-}\}$ corresponding to this choice of decoration.  Then $\mathcal{I}_{2}$ has two corresponding generators, which we will write $I_{+}$ and $I_{-}$, so $\mathcal{I}_{2} \cong \Z[q^{1/2},q^{-1/2}] I_{+} \oplus \Z[q^{1/2},q^{-1/2}] I_{-}$. \\
\ \\
\noindent For $\mathcal{I}_{4}$, there are twelve generators, depicted in Figure \ref{fig:I4Gens}. The inside disc $\lefty{D}_{L}$ is the shaded disk, while the outside disk $\righty{D}_L$ is the complement in the sphere. Thus, in terms of the combinatorial data the type $D$ generators have $\lefty{m}_{L} = (2,4)$ while $\righty{m}_{L} = (4,2)$. 

\myfig{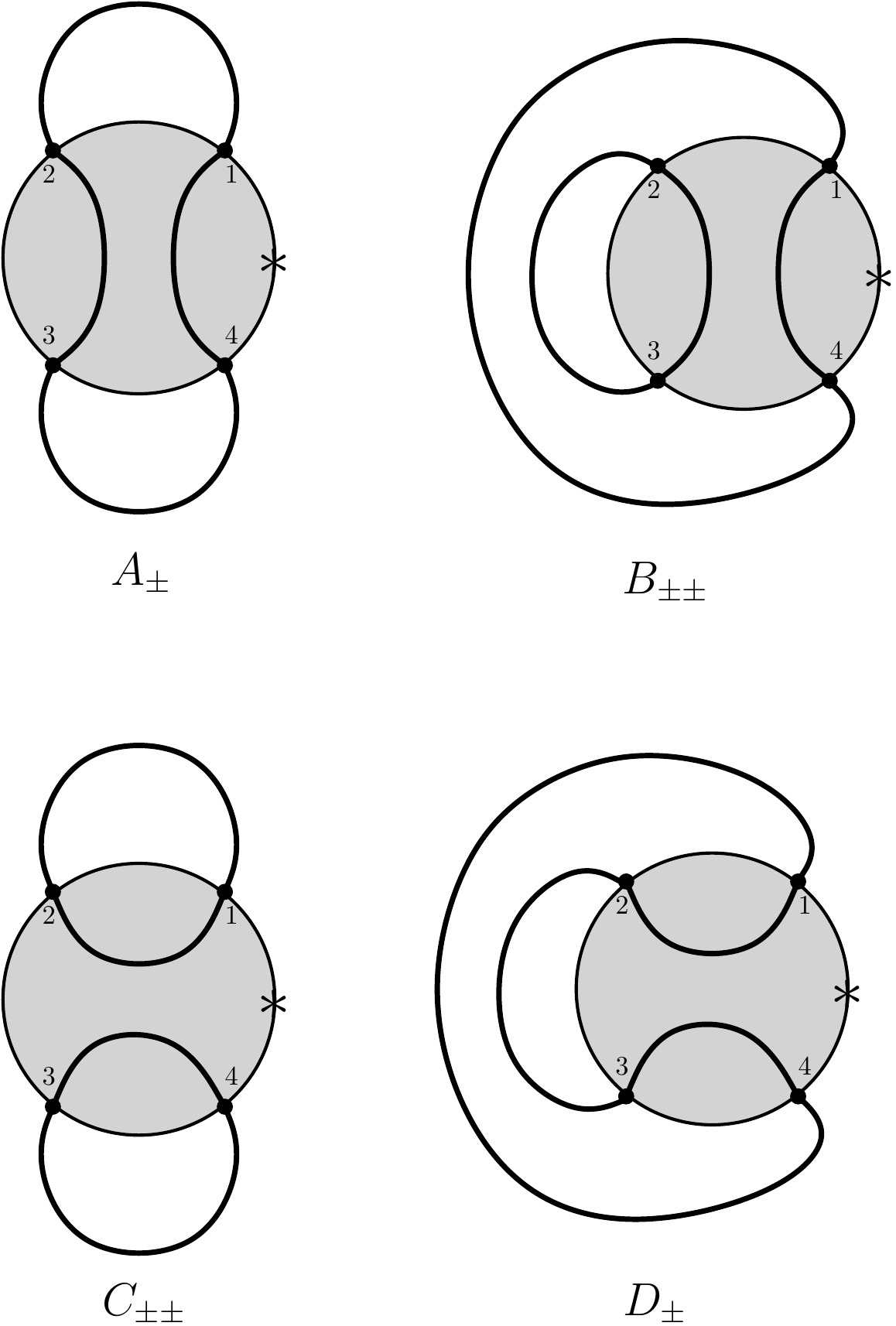}{The twelve generators of $\mathcal{I}_{4}$ grouped based on their inside and outside matchings. There a two generators of type $A$ and $D$, and four of type $B$ and $C$, as determined by the choice of decorations.}  

\noindent $\mathcal{I}_{6}$ has $104$ generators, which will not be listed here.

\section{Computing the decategorification of $\leftComplex{T}$}\label{sec:tangle}

\noindent Let $T$ be an oriented tangle diagram in an oriented two-dimensional disk $D_{T}$ with a marked point $\ast_{T} \in \partial D_{T}$, and the boundary circle oriented as the boundary of $D_{T}$. Let $\cross{T}$ be the set of crossings in $T$. Let $n_{\pm}(T)$ be the number of positive and negative crossings in $T$. In this section we explain how to compute $[\leftComplex{T}] \in K_{0}(\mathcal{B}\Gamma_{n})$ when $T$ has $2n$ endpoints. From the computation in section \ref{sec:decat} we know that
$$
[\leftComplex{T}] = \sum_{(L,\sigma) \in \cleaved{n}} [\leftComplex{T} I_{(L,\sigma)}] \in \mathcal{I}_{2n}
$$
We will describe the generators of $\leftComplex{T}$,  and the action of the idempotents, presently. We will also describe the bigradings for each generator. However, There are a few preliminaries before these descriptions. \\
\ \\
\noindent First, we will consider $D_{T}$ as being embedded in a sphere $S$, oriented compatibly with the orientation of $D_{T}$. Let $\righty{D}_{T} = S \backslash(\mathrm{int}\,D_{T})$. Along with the marked point $\ast_{T}$, this effects a decomposition of $S$ into inside and outside discs. 

\begin{defn}
An APS-resolution $\rho$ of $T$ is a map $\rho: \cross{T} \lra \{0,1\}$. For each APS-resolution, $\rho$, there is a planar diagram  $\rho(T)$ in $D^{2}$, called a APS-resolution diagram, where $\rho(T)$ is the diagram in $D^2$ obtained by locally replacing (disjoint) neighborhoods of the crossings of $T$ using the following rule for each crossing $c \in \cross{T}$:\\
$$
\inlinediag[0.3]{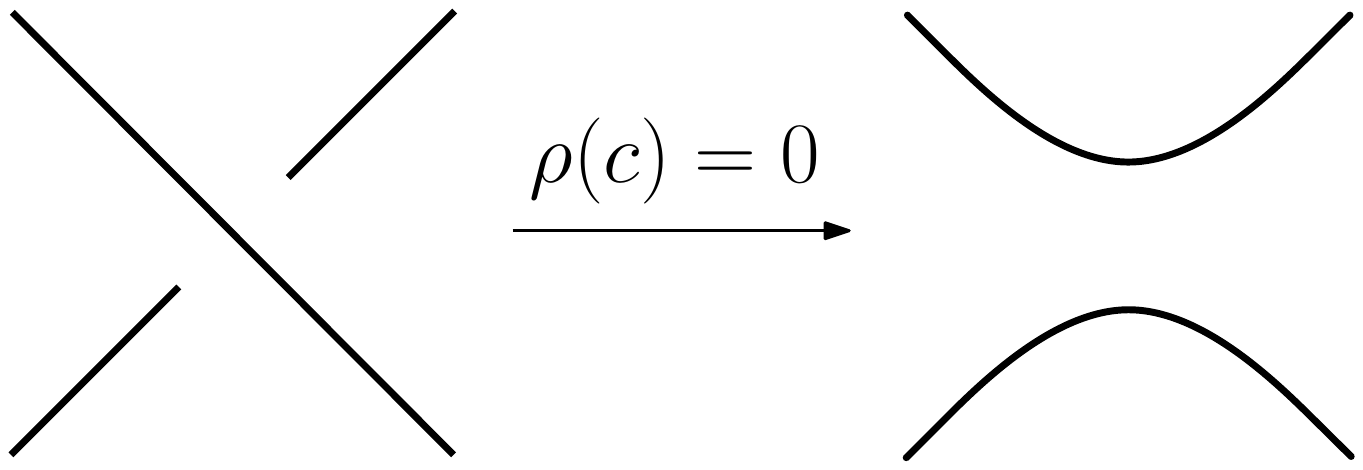} \hspace{.5in} \inlinediag[0.3]{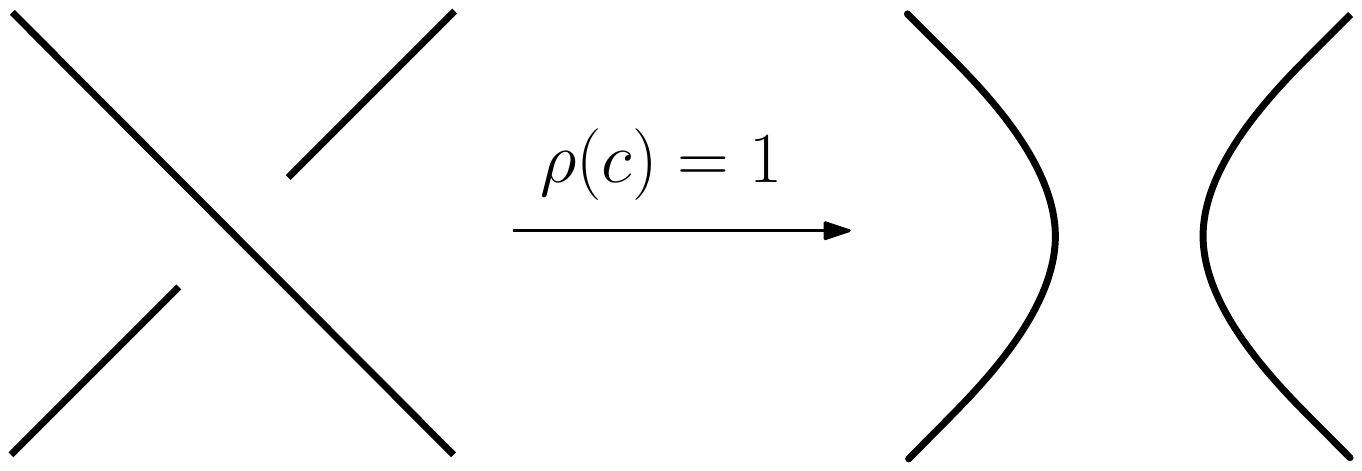} 
$$  
\end{defn}

\begin{defn}
A resolution of $T$ is a pair $(\rho, \righty{m})$ where $\rho$ is an APS-resolution of $T \subset D_{T}$, and $\righty{m}$ is a planar matching of $P_{T} = T \cap \partial D_{T}$ embedded in $\righty{D}_{T}$. The resolution diagram $\rho(T, \righty{m})$ for $(\rho, \righty{m})$ is the diagram in $S_{T}$ found by gluing $\rho(T)$ to $\righty{m}$ along $P_{T}$. The set of resolutions will be denoted $\resolution{T}$.
\end{defn}

\begin{defn}
For each $(\rho, \righty{m}) \in \resolution{T}$ the homological grading is $$h(\rho) = \sum_{c \in \cross{T}} \rho(c) - n_{-}(T)$$
\end{defn}

\begin{defn}
A {\em decorated} resolution for $T$ is a resolution $(\rho, \righty{m})$ and a map $s: \circles{\rho(T, \righty{m})} \lra \{+, -\}$
\end{defn}

\noindent The circles in $\rho(T, \righty{m})$ are of two types: 1) {\em free} circles -- those which do not intersect $P_{T}$, and 2) {\em cut} circles -- those which do. The set of free circles will be denoted $\freecircles{\rho(T, \righty{m})}$ while the cut circles will be denoted $\cutcircles{\rho(T, \righty{m})}$. 

\begin{defn}
The {\em quantum grading} of a decorated resolution $(\rho, \righty{m}, s)$ is
$$
I(\rho, \righty{m}, s) = h(\rho) + \sum_{C \in \freecircles{\rho(T, \righty{m})}} s(C) + \frac{1}{2}\sum_{C \in \cutcircles{\rho(T, \righty{m})}} s(C) + (n_{+}(T) - n_{-}(T))
$$
\end{defn}

\noindent The bigraded $\Z$-module obtained from $\leftComplex{T}$ by forgetting the differential, and the module structure over $\mathcal{B}\Gamma_{n}$, is generated by the decorated resolutions, shifted by the homological and quantum gradings:
$$
\leftComplex{T} = \bigoplus_{(\rho, \righty{m},s)} \Z\{(h(\rho), I(\rho, \righty{m},s)\}
$$

\noindent To describe the action of an idempotent $I_{(L,\sigma)}$ first notice that there is a map taking a decorated resolution to its {\em boundary} cleaved circle: $\partial(\rho, \righty{m}, s)$ is the cleaved circle $(\cutcircles{\rho(T, \righty{m})}, s|_{\mathrm{cut}})$ where $s|_{\mathrm{cut}}$ is the restriction of $s$ to the cut circles. This map is specified by the choice of marked point $\ast_{T}$; different choices will result in different maps.\\
\ \\
\noindent Then the action of $I_{(L,\sigma)}$ on $\leftComplex{T}$ is the linear extension of
$$
(\rho, \righty{m}, s) \cdot I_{(L,\sigma)} = \left\{\begin{array}{cl} (\rho, \righty{m}, s) & \partial(\rho, \righty{m}, s) = (L,\sigma) \\ 0 & \mathrm{Otherwise} \end{array} \right.
$$
Consequently, $\leftComplex{T} I_{(L,\sigma)}$ is spanned by the generators whose boundary cleaved circle is exactly $(L,\sigma)$. \\
\ \\
\noindent Since $\leftComplex{T} I_{(L,\sigma)}$ is a bigraded chain complex over $\Z$, we know how to compute $[\leftComplex{T}I_{(L,\Sigma)}]$. Let $P_{(L,\Sigma)}$ be the simple module over $\mathcal{B}\Gamma_{n}$ with $\Z$ in bigrading $(0,0)$ and such that $I_{(L,\sigma)}$ acts by the identity, while the action of the rest of $\mathcal{B}\Gamma_{n}$ is trivial. Each generator $(\rho, \righty{m}, s)$ gives a copy of $P_{(L,\Sigma)}$ in bigrading $(h(\rho), I(\rho, \righty{m}, s))$. Thus, 
$$
[\leftComplex{T}I_{(L,\Sigma)}] = \left(\sum_{\partial(\rho, \righty{m}, s) = (L,\sigma)} (-1)^{h(\rho)}q^{I(\rho, \righty{m}, s)} \right)P_{(L, \sigma)}
$$
Under the isomorphism with $\mathcal{I}_{2n}$, this implies that, in $K_{0}(\mathcal{B}\Gamma_{n})$,
$$
[\leftComplex{T}] = \sum_{(\rho, \righty{m}, s)}  (-1)^{h(\rho)}q^{I(\rho, \righty{m}, s)} I_{\partial (\rho, \righty{m}, s)}
$$
We illustrate with a few example computations.\\
\ \\
\noindent{\bf Example 1:} Suppose $T$ is a tangle in $D_{T}$ with no boundary, $P_{T} = \emptyset$. Then $\leftComplex{T}$ is the regular Khovanov homology, and $[\leftComplex{T}]$ is just $J_{T}(q)$. \\
\ \\
\noindent{\bf Example 2:} Suppose $P_{T} = \{p_{1}, p_{2}\}$ consists of two points, oriented as described in the section on cleaved links. Then the Grothendieck group is $\mathcal{I}_{2}$, which is two dimensional over $\Z[q^{1/2}, q^{-1/2}]$ with basis $I_{+}$ and $I_{-}$. There is only one right matching which can be used, and resolutions $(\rho, \righty{m}, s)$ come in pairs $r_{\pm}$ depending upon whether the single cut circle is adorned with a $+$ or $-$. Now $h(r_{+}) = h(r_{-})$ and there is a number $\widetilde{F}$ such that $I(r_{+}) = \widetilde{F} + 1/2$ while $I(r_{-}) = \widetilde{F} - 1/2$. Consequently, there is a polynomial $F(q)$ for which 
$$
[\leftComplex{T}] = F(q)\big(q^{1/2}I_{+} + q^{-1/2}I_{-}\big)
$$

\begin{center}
\begin{figure}
\includegraphics[scale=1]{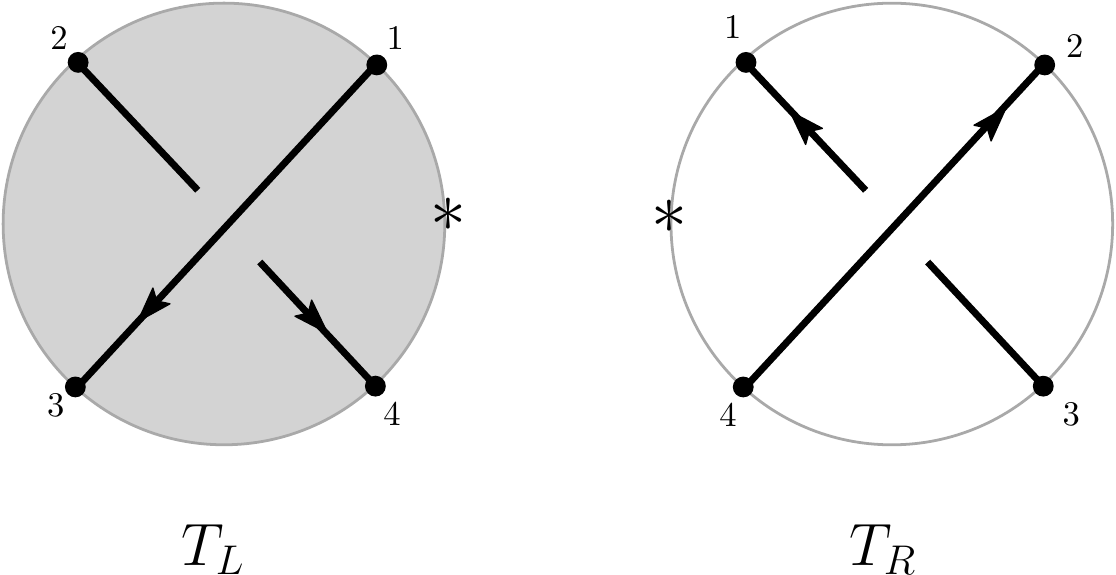}
\caption{The tangle $T_{L}$ on the left is an inside tangle, since the ordering of its endpoints is with the orientation on the boundary of the disk. It has a single positive crossing. The tangle on the right, $T_{R}$, also has a single positive crossing. However, it is an example of an outside tangle.}\label{fig:ex}
\end{figure}
\end{center}

\noindent{\bf Example 3:} We consider the tangle $T_{L}$ in Figure \ref{fig:ex} which has two arcs and a single positive crossing, with the choice of marked point as illustrated. The $0$ resolution consists of two vertical arcs, which can be glued to either outside matching $(24)$ or $(42)$. For $(0, (24), s)$ there is a single circle $C$ and thus two possibilities for $s$. We then have $h(0, (24), \pm) = 0 - 0 = 0$, since $n_{-} = 0$, and $I(0, (24), \pm ) = 0 + 0 \pm 1/2 + 1 - 0$. Furthermore, $\partial(0, (24), \pm ) =  A_{\pm}$ in the notation of Figure \ref{fig:I4Gens}. Thus, these generators contribute terms $q^{3/2}I_{A_{+}} + q^{1/2}I_{A_{-}}$. For the matching $(42)$ there are $2$ cut circles  $B_{\pm \pm}$. Once again, the homological grading is $0$, and the quantum grading is $I(0, (42),B_{\pm \pm}) = 0 + 0 \pm 1/2 \pm 1/2 + 1 - 0$, so we have terms $q^{2}I_{B_{++}} + q^{1}(I_{B_{+-}} + I_{B_{-+}}) + I_{B_{--}}$. For the $1$ resolution, we have two circles for matching $(24)$ which result in generators  $C_{\pm\pm}$. Now the homological grading is $1$ and the quantum grading is $I(1, (24),C_{\pm \pm}) = 1 + 0 \pm 1/2 \pm 1/2 + 1 - 0$. Therefore, we obtain terms $-q^{3}I_{C_{++}} - q^{2}(I_{C_{+-}} + I_{C_{-+}}) - qI_{C_{--}}$. Finally, for the $(42)$ matching we have $h(1,(42), D_{\pm}) = 1$ and $I(1, (42), D_{\pm}) = 1 + 0 \pm 1/2 + 1 - 0$, so we have terms $-q^{5/2}I_{D_{+}} - q^{3/2}I_{D_{-}}$. In total,
\begin{equation*}
\begin{split}
[\leftComplex{T_{L}}] &= q^{3/2}I_{A_{+}} + q^{1/2}I_{A_{-}} + q^{2}I_{B_{++}} + q^{1}(I_{B_{+-}} + I_{B_{-+}})\\
& + I_{B_{--}} - q^{3}I_{C_{++}} - q^{2}(I_{C_{+-}} + I_{C_{-+}}) - qI_{C_{--}} \\ 
& - q^{5/2}I_{D_{+}} - q^{3/2}I_{D_{-}}\\
\end{split}
\end{equation*}

\noindent Since the homotopy type of $\leftComplex{T}$ is an invariant of the tangle $T$, the decategorifications in $\mathcal{I}_{2n}$ are also invariants of the tangles. This can be proven directly, and the constructions above generalized, a process which we will return to in the sequel. 

\section{The decategorification of type D structures}\label{sec:typeD}

\noindent Above we mentioned that there are two type of invariants associated to a tangle. More specifically, suppose we have an an oriented tangle diagram $T$ in a oriented disc $D_{T}$ with a marked point $\ast_{T} \in \partial D_{T}$. Let $E = \mir{\partial D_{T}}$ be the boundary of $D_{T}$ with the {\em opposite} orientation. Let $P_{T} = \partial D_{T} \cap T$, ordered according to the orientation on $E$. Once again we think of $D_{T}$ as embedded in an oriented sphere $S$, but now $\lefty{D}_{T} = S \backslash(\mathrm{int}\,D_{T})$ is the inside disc for this decomposition of $S$. \\
\ \\
\noindent To this configuration is associated a different type of algebraic invariant: a type $D$ structure $\rightComplex{T}$. This is a $(1,0)$ map $\delta: \rightComplex{T_{2}} \lra  \mathcal{B}\Gamma_{n} \otimes_{I}  \rightComplex{T}$ which satisfies a certain structural identity that will not play a role in this paper. Type D structures also have Grothendieck groups, and, following Ina Petkova in \cite{Petk} and the results of Lipshitz, Oszv\'ath, and Thurston in  \cite{Mor}, the computations are essentially the same as above\footnote{The various bimodule constructions of \cite{Mor} showing the equivalence of categories can be adapted to the bordered Khovanov setting}. Thus, the previous section more or less tells us how to compute the decategorification of the type $D$ structure. \\
\ \\
\noindent However, the point of introducing the two structures was to obtain a gluing $\leftComplex{T_{1}} \boxtimes \rightComplex{T_{2}}$, so we will adapt the definition of the decategorification of the type $D$ structure to reflect this pairing in the decategorifications. \\
\ \\
\noindent We will think of the decategorification of $\leftComplex{T_{1}}$ as the map $\Z[q^{1/2}, q^{-1/2}] \lra K_{0}(\mathcal{B}\Gamma_{n}) \cong \mathcal{I}_{2n}$ defined by $f(q) \lra f(q)[\leftComplex{T_{1}}]$, where $|T_{1} \cap \partial D_{T_{1}}| = 2n$. Accordingly, we will think of the decategorification of $\rightComplex{T_{2}}$ as an element in the dual $\mathcal{I}_{2n} \cong K_{0}(\mathcal{B}\Gamma_{n}) \lra \Z[q^{1/2},q^{-1/2}]$, and thus as an element $[\rightComplex{T_2}]$ in  $\mathcal{I}_{2n}^{\ast} = \textsc{Hom}_{\Z[q^{\pm 1/2}]}(\mathcal{I}_{2n}, \Z[q^{\pm 1/2}])$. \\
\ \\
\noindent We will describe this map somewhat tersely: the generators of $\rightComplex{T_{2}}$ are the triples $(\rho, \lefty{m}, s)$ where $\rho$ is an APS-resolution of the diagram $T_{2}$, $\lefty{m}$ is a planar matching of $P_{T_{2}}$ in $\lefty{D}_{T_{2}}$ (an {\em inside} matching) and $s: \circles{\rho(\lefty{m}\#T_{2})} \lra \{+,-\}$. The homological and quantum gradings of this generator are computed identically to those in $\leftComplex{T_{1}}$, described above. Likewise, each generator has a boundary $\partial(\rho, \lefty{m}, s) \in \cleaved{n}$ obtained by erasing all the free circles. The action of the idempotent $I_{(L,\Sigma)}$ on the {\em left} of $\rightComplex{T_{2}}$ is trivial on generators whose boundary is different from $(L,\sigma)$ and the identity on those generators for which $(L,\sigma)$ is the boundary. \\
\ \\
\noindent The map $[\rightComplex{T_2}] \in \textsc{Hom}_{\Z[q^{\pm 1/2}]}(\mathcal{I}_{2n}, \Z[q^{\pm 1/2}])$ is the linear extension of the following map on generators $I_{(L,\sigma)}$ of $\mathcal{I}_{2n}$:
$$
[\rightComplex{T_2}](I_{(L,\sigma)}) = \sum_{\partial(\rho,\lefty{m},s)= (L,\sigma)} (-1)^{h(\rho)}q^{I(\rho,\lefty{m},s)}
$$

\noindent{\bf Example:} We consider the tangle $T_{R}$ in Figure \ref{fig:ex} which is an outside disc containing a single positive crossing. To specify the map $[\rightComplex{T_{R}}]$ we need to specify the image on each of the generators in Figure \ref{fig:I4Gens}.  Note, however, that the arcs coming from resolving the crossing will be those {\em outside} the shaded discs in Figure \ref{fig:I4Gens}. Thus a $0$-resolution of the crossing gives $(42)$ as the {\em outside} matching. If we choose $(24)$ for the inside matching, then we obtain the generators $D_{\pm}$. For $D_{+}$ there is one positive circle, so $I(D_{+}) = 0 + 1/2 + 1$ and $h(D_{+}) = 0 - 0$. Thus $[\rightComplex{T_{R}}](D_{+}) = q^{1/2+1} = q^{3/2}$. However, if we choose $(42)$ as the inside matching, then we obtain the generators $B_{\pm\pm}$. Since $B_{-+}$ has $h(B_{-+}) = 0 - 0$ and $I(B_{-+}) = 0 - 1/2 + 1/2 + 1 = 1$ we can compute $[\rightComplex{T_{R}}](B_{-+}) = q$. To obtain $C_{++}$ we need to choose $(24)$ as both the inside and outside matching. Thus we need the $1$ resolution of the tangle diagram. Then $h(C_{++}) = 1$ and $I(C_{++}) = 1 + 1/2 + 1/2 + 1 = 3$.  Similar computations for each of the generators produce the following map:

$$
\begin{array}{lclcl}
A_{+} \lra -q^{5/2} & \hspace{0.75in} & B_{-+} \lra q & \hspace{0.75in} & C_{-+} \lra -q^2 \\
A_{-} \lra -q^{3/2} & \hspace{0.75in} & B_{--} \lra 1 & \hspace{0.75in} & C_{--} \lra -q\\
B_{++} \lra q^2 & \hspace{0.75in} & C_{++} \lra -q^3& \hspace{0.75in} & D_{+} \lra q^{3/2}\\
B_{+-} \lra q & \hspace{0.75in} & C_{+-} \lra -q^2 & \hspace{0.75in} & D_{-} \lra q^{1/2} \\
\end{array}
$$

\noindent Note that the map $[\rightComplex{T}]$ is in principal determined by the image of the cleaved links with all $+$ decorations, since changing a $+$ to a $-$ does not change the APS-resolutions and matching for those states. It does multiply by $q^{-1}$ due to the change in quantum grading. Thus, once we know that $B_{++} \lra q^{2}$ we know that $B_{-+} \lra q^{-1}(q^{2})$, for example. 

\section{Recovering the Jones polynomial}\label{sec:Jones}

\noindent Suppose $S$ is an oriented sphere decomposed as $\lefty{D} \cup \righty{D}$, where $\lefty{D}$ and $\righty{D}$ are two oriented discs who intersect only on their common boundary, and that boundary is oriented as the boundary of $\lefty{D}$. Let $\ast \in \partial \lefty{D}$. Suppose further that $L$ is an oriented link diagram in $S$ which intersects $\partial \lefty{D}$ transversely away from its crossings. Then $T_{1} = L \cap \lefty{D}$ and $T_{2} = L \cap \righty{D}$ are an inside tangle and an outside tangle, respectively. 

\begin{prop}
The composition $$\Z[q^{1/2},q^{-1/2}] \stackrel{[\leftComplex{T_{1}}]}{\lra} \mathcal{I}_{2n} \stackrel{[\rightComplex{T_{2}}]}{\lra} \Z[q^{1/2},q^{-1/2}]$$ equals multiplication by the Jones polynomial $J_{L}(q)$, described in the introduction.
\end{prop}

\begin{center}
\begin{figure}
\includegraphics[scale=0.75]{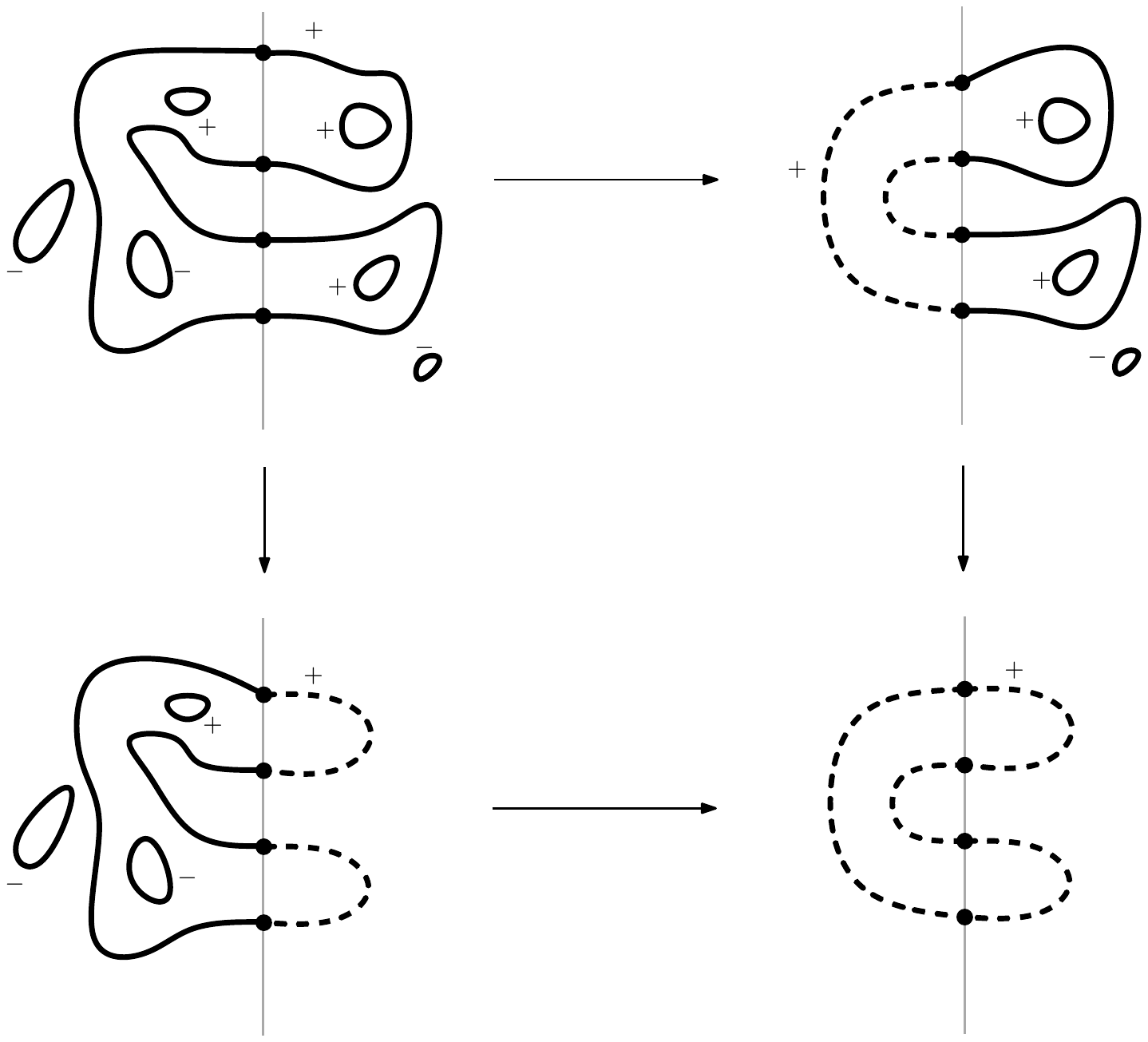}
\caption{The top left diagram represents a decorated resolution $\xi$ used to compute $J_{L}(q)$. The bottom left is the corresponding generator in $\leftComplex{T_{1}}$. The dashed arcs in the bottom left represent the outside matching. They are the arcs which remain when we delete all the free circles to the right of the vertical line. The top right is the generator  used in $\rightComplex{T_{2}}$.The bottom right is the cleaved link that is the common boundary of these generators. Notice that the bottom right diagram is determined solely by the cut circles, and their decorations, shown in the top left diagram. Notice also that the top right and bottom left can be identified along their boundary to recover the diagram in the top left, simply by laying one diagram atop the other so that the boundaries match.}\label{fig:Jones}
\end{figure}
\end{center}

\noindent{\bf Proof:} $J_{L}(q)$ is the sum $\sum_{(\rho,s)} (-1)^{h(\rho)} q^{I_{J}(\rho,s)}$ where  1) $\rho: \cross{L} \lra \{0,1\}$ is a resolution of the diagram $L$ with diagram $\rho(L)$ found using same rules as for APS-resolutions, and 2) $s: \circles{\rho(L)} \lra \{+,-\}$ as above. The value $h(\rho)$ is computed as above, using all the crossings in $L$. $I_{J}(\rho,s)$ is computed as above, noting that there will be no cut circles in $\rho(L)$ so every circle will contribute $\pm 1$ to the count. \\
\ \\
\noindent For each such resolution diagram $\rho(L)$, let $\lefty{\rho}(L) = \rho(L) \cap \lefty{D}$. Then $\lefty{\rho}(L)$ is an APS-resolution diagram for $T_{1}$, whose resolution is simply the restriction of $\rho$ to the crossings in $\lefty{D}$. A similar story holds for $\rho(L) \cap \righty{D} = \righty{\rho}(L)$. Let $\righty{m}$ be the planar matching obtained by deleting the free circles from $\righty{\rho}(L)$ and let $\lefty{s}$ be the decoration on $\lefty{\rho} \# \righty{m}$ found by restricting $s$. Then $(\lefty{\rho}, \righty{m}, \lefty{s})$ is a generator for $\leftComplex{T_{1}}$. If we take $\lefty{m}$ to be the planar matching found by deleting the free circles of $\lefty{\rho}(L)$ then $(\righty{\rho}, \lefty{m}, \righty{s})$ is a generator of $\rightComplex{T_{2}}$. \\
\ \\
\noindent Thus every generator $\xi$ of $\complex{L}$ maps to a pair $(\lefty{\xi}, \righty{\xi})$ of generators from $\leftComplex{T_{1}}$ and $\rightComplex{T_{2}}$. Furthermore, $\partial \lefty{\xi} =$ $\partial \righty{\xi} =$ $(\lefty{m}\#\righty{m}, \sigma)$ where $\sigma$ is the restriction of $s$ to the circles of $\rho(L)$ intersecting $\partial \lefty{D}$. This is illustrated in Figure \ref{fig:Jones}. In fact, given any such pair of generators $(\lefty{\xi}, \righty{\xi})$ with $\partial \lefty{\xi} = \partial \righty{\xi}$ there is a unique generator $\xi$ of $\complex{L}$ which maps to the pair.  Thus, this construction generates a bijection between the generators of $\complex{L}$ and the set of pairs $(\lefty{\xi}, \righty{\xi})$ with $\partial \lefty{\xi} = \partial \righty{\xi}$.
\noindent It is straightforward to check that, under this bijection, $h(\rho) = h(\lefty{\rho}) + h(\righty{\rho})$, since the crossings of $L$ are partitioned between $\lefty{D}$ and $\righty{D}$. It is also true that $I_{J}(\rho,s)= I(\lefty{\rho}, \righty{m}, \lefty{s}) + I(\righty{\rho}, \lefty{m}, \righty{s})$ since a $\pm$ cut circle will be counted as $\pm 1/2$ in the sum defining each of $I(\lefty{\rho}, \righty{m}, \lefty{s})$ and $I(\righty{\rho}, \lefty{m}, \righty{s})$, but as $\pm 1$ in the sum for $I_{J}(\rho,s)$. Thus the monomial in $J_{L}(q)$ for each generator $\xi$ of $\complex{L}$ is a product
$$
(-1)^{h(\xi)} q^{I_{J}(\xi)} = (-1)^{h(\lefty{\xi})} q^{I(\lefty{\xi})}(-1)^{h(\righty{\xi})} q^{I(\righty{\xi})} 
$$
where the factors on the right are monomial terms in $[\leftComplex{T_{1}}]$ and $[\rightComplex{T_{2}}]$. Furthermore, provided that $\partial \lefty{\xi} = \partial \righty{\xi}$ the product of the terms on the right will equal a monomial in $J_{L}(q)$. \\
\ \\
\noindent Let $I^{\ast}_{(L,\sigma)}$ be the Kronecker dual functional for $I_{(L,\sigma)}$. Then the composition 
$$\Z[q^{1/2},q^{-1/2}] \stackrel{[\leftComplex{T_{1}}]}{\lra} \mathcal{I}_{2n} \stackrel{[\rightComplex{T_{2}}]}{\lra} \Z[q^{1/2},q^{-1/2}]$$
is multiplication by the polynomial that results from
$$
\left(\sum_{I_{(L,\sigma)}}\left(\sum_{\partial(\righty{\rho}, \lefty{m}, s')=(L,\sigma)} (-1)^{h(\righty{\rho})} q^{I(\righty{\rho}, \lefty{m}, s')}\right) I^{\ast}_{(L,\sigma)}\right) \circ \left(\sum_{I_{(L,\sigma)}}\left(\sum_{\partial(\lefty{\rho}, \righty{m}, s)=(L,\sigma)} (-1)^{h(\lefty{\rho})} q^{I(\lefty{\rho}, \righty{m}, s)}\right) I_{(L,\sigma)}\right) 
$$
$$
= \sum_{\partial\lefty{\xi} = \partial\righty{\xi}} (-1)^{h(\lefty{\xi})} q^{I(\lefty{\xi})}(-1)^{h(\righty{\xi})} q^{I(\righty{\xi})} 
$$
where $\lefty{\xi}$ and $\righty{\xi}$ are allowed to range over all generators of $\leftComplex{T_{1}}$ and $\rightComplex{T_{2}}$. \\
\ \\
\noindent From our argument above, the terms in this sum are in one-to-one correspondence with the terms in the sum defining $J_{L}(q)$, and equal the corresponding term. Thus, the sum is equal to $J_{L}(q)$. $\Diamond$\\
\ \\
\noindent{\bf Example:} We can glue the discs containing $T_{L}$ and $T_{R}$ along their boundaries so that the marked and labeled points are identified. This gives a link diagram for the positive Hopf link. Computing $\rightComplex{T_{R}} \circ \leftComplex{T_{L}}$ gives
\begin{equation*}
\begin{split}
\rightComplex{T_{R}}(q^{3/2}I_{A_{+}}& + q^{1/2}I_{A_{-}} + q^{2}I_{B_{++}} + q(I_{B_{+-}} + I_{B_{-+}})+ I_{B_{--}}\\ 
& - q^{3}I_{C_{++}} - q^{2}(I_{C_{+-}} - I_{C_{-+}}) - qI_{C_{--}} - q^{5/2}I_{D_{+}} - q^{3/2}I_{D_{-}})\\
& = q^{3/2}(-q^{5/2}) + q^{1/2}(-q^{3/2}) + q^{2}(q^2) + q(q + q)\\
& +1 - q^{3}(-q^3) - q^{2}(-q^2 - q^2) - q(-q)  - q^{5/2}(q^{3/2}) - q^{3/2}(q^{1/2})\\
& = -q^4 - q^2 + q^4 + 2q^2 + 1 + q^6 + 2q^4 + q^2 - q^4 - q^2\\
& = q^{6} + q^{4} + q^{2} + 1 
\end{split}
\end{equation*}
This is the correct polynomial $J_{L}(q)$ for this Hopf link.

\begin{center}
\begin{figure}
\includegraphics[scale=0.6]{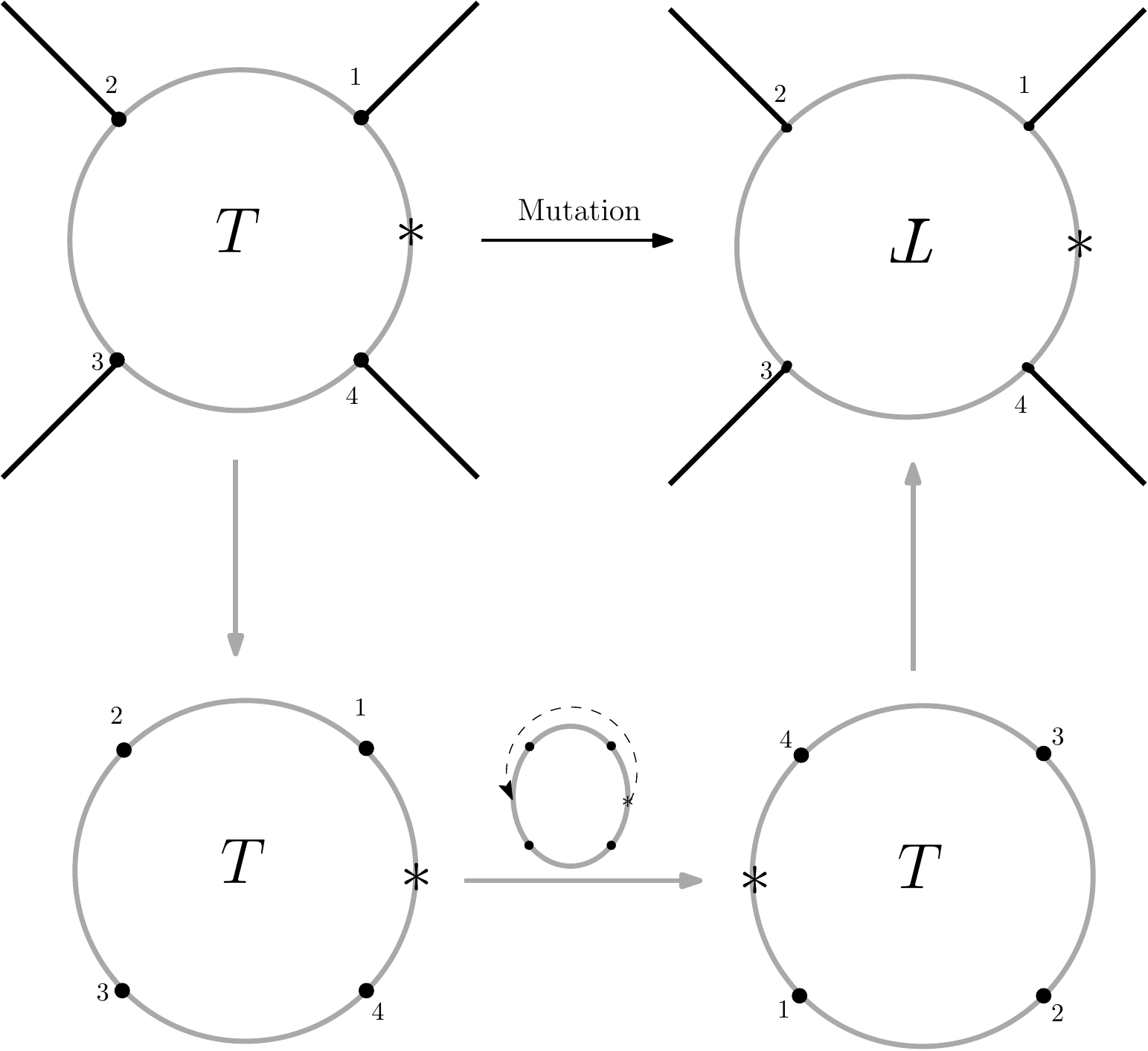}
\caption{The top row is a mutation of the diagram, where we have only depicted the region around the tangle $T$. We have chosen a marked point and labeled the intersections with the boundary of the local disc, so that the local disc is an inside disc. The bottom left is the inside disc. Moving the marked point two segments counter-clockwise is shown in the lower row. This relabels the intersection points as well. Gluing the new inside disc according to the boundary data produces the mutant.}\label{fig:mutation}
\end{figure}
\end{center}

\section{Mutation Invariance}\label{sec:mutation}

\noindent We use the results of the last section to provide an alternative proof of the well known fact that $J_{L}(q)$ is invariant under mutations of $L$. A diagram $L'$ is a mutant of a diagram $L$ if there is a disc $D$ in $L$ which intersects $L$ in $4$ points, such that $L'$ is obtained by removing $D$ and gluing it back with a $180^{\circ}$ twist about one of the three principal axes in $\R^{3}$. Using the constructions in this paper, we can describe mutation in terms of moving the marked point $\ast$ used in gluing the local disc. Pick $\ast$ on $\partial D$ and let $T_{1} = L \cap D$ and $T_{2} = L \cap (S \backslash D)$. With the marked point, $T_{1}$ is an inside tangle in $D$, and $T_{2}$ is an outside tangle. To obtain $L'$ we move $\ast$ on $\partial D$ two segments counter-clockwise, while keeping the marked point on $\partial(S \backslash D)$ fixed. This process, and how it achieves mutation are shown in Figure \ref{fig:mutation}. Gluing to align the marked points results in the mutant diagram $L'$. \\
\ \\
\noindent In \cite{Wehr}, it is shown that this one type of mutation is enough to obtain mutations where we rotate the diagram in $D$ around the other two perpendicular axes. These rotations do not have easy interpretations in terms of the marked point. \\
\ \\
\noindent Consider the element $[\leftComplex{T_{1}}] \in \mathcal{I}_{4}$. Given a generator $(\rho, \righty{m}, s)$ of $\leftComplex{T_{1}}$, moving the marked point $\ast$ changes $\partial(\rho, \righty{m}, s)$ but not $h(\rho)$, hence we will reverse orientations on all the strands and thus not change any of the crossings, or $I(\rho, \righty{m}, s)$. Consequently, we can think of moving $\ast$ as providing a map $M_{\ast}: \mathcal{I}_{4} \lra \mathcal{I}_{4}$ under which $[\leftComplex{T_{\ast}}]$ is $M_{\ast}[\leftComplex{T}]$ if $T_{\ast}$ is the inside tangle with the new marked point. From the generators in Figure \ref{fig:I4Gens} we can compute the effect on the generators of $\mathcal{I}_{4}$ of moving $\ast$ {\em one} segment counter-clockwise:

$$
\begin{array}{lclclcl}
A_{\pm} \lra D_{\pm} & \hspace{0.3in} & D_{\pm} \lra A_{\pm} & \hspace{0.3in} & B_{++} \lra C_{++} & \hspace{0.3in} & C_{++} \lra B_{++} \\
B_{--} \lra C_{--} & \ & C_{--} \lra B_{--} & \ & C_{+-} \lra B_{+-} & \ & C_{-+} \lra B_{-+} \\
B_{+-} \lra C_{-+} & \ & B_{-+} \lra C_{+-} & \ & \ & \ & \ \\
\end{array}
$$

\noindent Notice that the signs in $B_{+-} \lra C_{-+}$ are switched. Recall that the circles in $B_{+-}$ are ordered by the lowest number appearing on the circle. This number can change when we move $\ast$, but the decoration on the circle does not. \\
\ \\
\noindent To compute $M_{\ast}$ we must move $\ast$ two segments. This produces a map on $\mathcal{I}_{4}$ which fixes all the generators except $B_{+-} \leftrightarrow B_{-+}$ and $C_{+-} \leftrightarrow C_{-+}$. Thus, to prove mutation invariance it suffices to see that the coefficient of $B_{+-}$ in  
$[\leftComplex{T_{1}}]$ is the same as $B_{-+}$, and likewise for $C_{+-}$ and $C_{-+}$. For then, when we apply $\rightComplex{T_{2}}$ to $\leftComplex{T_{1}}$ we will get the same result as applying it to $M_{\ast}^{2}[\leftComplex{T_{1}}]$. As these compositions compute the Jones polynomial for $L$ and its mutant $L'$, we will have verified the mutation invariance. \\
\ \\
\noindent Suppose we have a generator $(\rho, \righty{m}, s)$ of $\leftComplex{T_{1}}$ with $\partial(\rho, \righty{m}, s) = B_{-+}$. Then we can find a corresponding generator $(\rho', \righty{m}', s')$ whose boundary is $B_{+-}$. We take $\rho'=\rho$ and $\righty{m}'=\righty{m}$, and thus have the same resolution diagram. All that changes is that $s'$ will evaluate on the cut circles opposite how $s$ does. Therefore, $h(\rho) = h(\rho')$ and $I(\rho, \righty{m}, s) = I + 1/2 - 1/2$, where we have pulled out the contribution of the cut circles, while $I(\rho', \righty{m}', s') = I - 1/2 + 1/2$. As these are equal, the coefficient on $B_{-+}$ from $(\rho, \righty{m}, s)$ is identical to that on $B_{+-}$ from $(\rho', \righty{m}', s')$. Since the map $(\rho, \righty{m}, s) \rightarrow (\rho', \righty{m}', s')$ is evidently a bijection from the generators with boundary $B_{-+}$ to those with boundary $B_{+-}$, we can conclude that the coefficient on $B_{-+}$ in $[\leftComplex{T_{1}}]$ is identical to that on $B_{+-}$. A similar argument shows that the coefficients on $C_{+-}$ and $C_{-+}$ are also equal. $\Diamond$

\end{document}